\documentclass[11pt]{amsart}
\usepackage[paperwidth=7in, paperheight=10in, margin=.875in]{geometry}
 \usepackage[backref,colorlinks,linkcolor=red,anchorcolor=green,citecolor=blue]{hyperref}
\usepackage{amsfonts,amssymb,bm}
\usepackage{amsmath,fancyvrb}
\usepackage{graphicx}
\usepackage{mathrsfs} 
\usepackage{graphicx}
\usepackage{cite}
\usepackage{enumerate}
\usepackage{tikz}
\usepackage{caption}
\usepackage{easyReview}
\numberwithin{equation}{section}
\newtheorem{thm}{Theorem}[section]

\newtheorem{example}[thm]{Example}
\newtheorem{defn}[thm]{Definition}

\DeclareMathOperator{\Tr}{Tr}

\begin{document}
 \title[]{On the generation of  necklaces and bracelets in R}

\author{Elvira Di Nardo$^\star$}
\address{$\star$Dipartimento di Matematica \lq\lq G. Peano\rq\rq, Università degli Studi di Torino, Via Carlo Alberto 10, 10123 Torino, Italy}
\email{elvira.dinardo@unito.it}
\urladdr{elviradinardo.it}  

\author{Giuseppe Guarino$^\dag$}
\address{$\dag$Medical School, Università Cattolica del Sacro Cuore, 
  Largo Agostino Gemelli 8, 00168, Roma, Italy}
\email{giuseppe.guarino@rete.basilicata.it}
 
\maketitle

\begin{abstract}
This note introduces a code snippet in {\tt R} aiming to generate necklaces as well as bracelets. Among various uses, necklaces are useful tools to manage traces of products of random matrices.  Functionality for necklaces and bracelets is provided with some examples of applications such as Lyndon words and de Bruijn sequences. The routines are collected in the {\tt Necklaces} package  available from the Comprehensive {\tt R} Archive Network.
\end{abstract}

 {\small {\bf keywords:} necklace, bracelet, multi-index composition, Lyndon word, de Bruijn sequence.}

\section{Introduction} \label{sec:intro}

Necklaces and bracelets are fundamental combinatorial objects arising in the field of text algorithms \cite{kociumaka2014computing}, in the construction of single-track Gray codes \cite{schwartz1999structure}, in the analysis of circular DNA and splicing systems \cite{de2017unavoidable} or in managing the distribution of individual wire faults in parallel on-chip links \cite{vitkovskiy2012combinatorial}. Applications in probability and statistics are related to random splitting \cite{alon2021random}, necklace processes  \cite{mallows2008necklace} or  fixed-density $m$-ary  necklaces \cite{ruskey1999efficient}.
Dealing with matrix algebra,  they are used for the well-known invariance property of the trace, that is the trace of a matrix product is invariant under cyclic permutations. For example if $\Sigma$ is a $p \times p$ non-negative definite matrix 
and $H_1, H_2 \in {\mathbb C}^{p \times p}$  are complex matrices, then
$\Tr[(\Sigma H_1)^2 (\Sigma H_2)] = \Tr[(\Sigma H_1) (\Sigma H_2) (\Sigma H_1) ] = \Tr[ (\Sigma H_2)(\Sigma H_1)^2].$ The set of indexes $\{112,121,211\}$ is an example of binary necklace of  length $3$ on the alphabet $\{1,2\}.$ 

More in details, a $m$-ary necklace of length $n$ is an equivalence class of $n$-character strings over an alphabet $X$ of size $m,$ taking all rotations  equivalent \cite{flajolet2009analytic}. That is, if $\sigma$ denotes the rotation pictured as $x_1 \mapsto x_2, x_2 \mapsto x_3, \ldots, x_n \mapsto x_1,$ the $m$-ary $n$-tuples $y$ and $z$ are in the same  necklace if and only if $\sigma^j(y)=z$ for some integer $j.$ Sometimes a $m$-ary necklace of length $n$ is represented with $n$ circularly connected beads coloured with up to $m$ different colours. Observe that the $n$-strings of a necklace are $n$-tuples of the set $X^n.$ Without losing generality, the elements of $X$ can be labelled with the first $m$ positive integers, that is $X = [m] = \{1,2,\ldots,m\},$ or with  $X = \{0,1,2,\ldots,m-1\}$ including the zero. 
A $m$-ary bracelet of length $n$ is an equivalence class where not only all rotations are equivalent but also all reflections (reverse order).   

It is quite natural to choose the lexicographically smallest string $\mathfrak{a}$ as the representative of a necklace (or bracelet).  Let us denote the class with $\langle \mathfrak{a}\rangle.$ 
\begin{example}\label{ex1}
{\rm There are $6$ binary necklaces of length $4:$
$$\langle 0\,0\,0\,1 \rangle=\{1000,0100,0010,0001\}, \langle 0\,0\,1\,1\rangle =\{1100,1001,0011,0110\},$$
$$\langle 0\,1\,0\,1\rangle=\{1010,0101\}, \langle 0\,1\,1\,1\rangle=\{1110,1101,1011,0111\},$$
$$\langle 1\,1\,1\,1\rangle=\{1111\},\langle 0\,0\,0\,0\rangle=\{0000\}, $$
all of which are also binary bracelets.   For $n = 6,$ we have 
$14$ binary necklaces but only $13$ binary bracelets as the necklaces obtained rotating $101100$ and  $001101$  are reversal of each other, resulting into a single bracket class, see Example \ref{3.4} in Section $3$.}
\end{example}
Necklaces turn to be useful also in free probability \cite{Speicher} when computing moments of suitable products of random matrices. For example, using necklaces  a closed form expression is given in \cite{MR3163834} to compute 
\begin{equation}
{\mathbb E}\bigg({\rm Tr}[W(n) H_1]^{i_1} \cdots {\rm Tr}[W(n) H_m]^{i_m}\bigg)
\label{genmom}
\end{equation}
where $H_1, \ldots, H_m \in {\mathbb C}^{p \times p}$ are complex matrices, $i_1, \ldots, i_m$ are non-negative integers and $W(n)$ is a non-central Wishart square random matrix of order $p$ and parameter $n.$ 
\begin{example}\label{ex0} {\rm For example\footnote{A {\tt Maple} procedure to compute these moments is available in \cite{DiNardo1}.} for $m=2$ and $i_1=1, i_2=1$ we have
\begin{eqnarray*}
&   & {\mathbb E}\big(\Tr[W(n) H_1] \Tr[W(n) H_2] \big) = \Tr(\Omega \Sigma H_1) \Tr(\Omega \Sigma H_2)  -  \Tr(\Omega \Sigma H_1 \Sigma H_2) \\
& - & \Tr(\Omega \Sigma H_2 \Sigma H_1) + n^2 \Tr(\Sigma H_1) \Tr(\Sigma H_2) + n
\Tr(\Sigma H_1 \Sigma H_2) - n \Tr(\Omega \Sigma H_1) \times \\
& \times &  \Tr(\Sigma H_2)  -  n \Tr(\Sigma H_1) \Tr(\Omega  \Sigma H_2),
\end{eqnarray*}
where $\Sigma$ is the covariance matrix and $\Omega$ is the non-centrality matrix of the non-central Wishart random matrix $W(n).$}
\end{example}
In computing \eqref{genmom}, the notion of {\sl type} of a necklace
has been considered, according to the occurrence of each character of the alphabet in the $n$-tuple. Suppose  $\bm i = (i_1, i_2,  \ldots,  i_m)$ and $X=[m].$
\begin{defn}\label{def1.1}
A necklace is said to be of type $\bm i$ if the character $1$ appear $i_1$ times in the $n$-tuples, $2$ appear $i_2$ times  and so on. 
\end{defn}
Therefore a necklace of type $\bm i$ is such that $i_j$ denotes the number of times the non-negative integer $j$ occurs in a string of the necklace. The choice  $X=[m]$ is not a constrain as $1$ can be replaced by whatever integer,  $2$ by the subsequent one and so on.  A similar definition can be given for bracelets. 

In {\tt R}, the only routines managing necklaces and bracelets are in the {\tt numbers} package \cite{numbers}, where the Polya’s enumeration theorem \cite{harary1994polya} is used to compute the 
number of representatives of  $m$-ary necklaces  (resp. bracelets)  of length $n,$ see Example \ref{3.4}. Instead the {\tt Necklaces} package can be used not only for producing the lists of all the $n$-tuples in a $m$-ary necklace (resp. bracelet)  of length $n,$ but also the lists of all the representatives  of $m$-ary necklaces (resp. bracelets) of type $\bm i$  as well as of $m$-ary necklaces  (resp. bracelets)  of length $n.$   

This note is organized as follows. Next section describes the functionality of those routines in the {\tt kStatistics} package called by the routines of the {\tt Necklaces} package.
Section $3$ adds more details  about  functionality of the routines of the {\tt Necklaces} package as well as the implemented methods. Their performance is illustrated with many examples. Section $4$ is devoted to examples of applications: the generation of Lyndon words and de Bruijn sequences. Some concluding remarks end the paper.  

%------------------------------------------------------------------------
\section{Functions of the kStatistics package} \label{sec:2} 
%------------------------------------------------------------------------
In order to use the functions of the {\tt Necklaces} package, first install the {\tt kStatistics}  package \cite{kStatistics}. Indeed the {\tt Necklaces} package calls the {\tt nPerm} function of the {\tt kStatistics}  package. This function returns all the permutations of an input vector of non-negative integers.
\begin{example}
{\rm  Run {\tt nPerm(c(1,0,1,1))}  to get
$\{1 1 1 0, 1 1 0 1, 1 0 1 1, 0 1 1 1\}.$}
\end{example}
\begin{example}
{\rm  Run {\tt nPerm(c(0,1,2))}  to get
$\{0 1 2, 0 2 1, 1 0 2, 1 2 0, 2 0 1, 2 1 0 \}.$}
\end{example}
The {\tt Necklaces} package calls also  the {\tt mkT} function, which returns all the compositions of a multi-index $\bm i = (i_1,  i_2,  \ldots,  i_m)$ in $n$ multi-indexes  \cite{di2022kstatistics}, that is all the lists $({\bm v}_1, \ldots, {\bm v}_n)$ of $m$-length multi-indexes such that ${\bm v}_1 + \cdots + {\bm v}_n = \bm i.$
\begin{example} {\rm To get all the compositions of the multi-index 
$(1,0,1)$ in two multi-indexes of length $3$, run}
\begin{verbatim}
> mkT(c(1,0,1),2, TRUE)
[( 0 0 1 )( 1 0 0 )]
[( 1 0 0 )( 0 0 1 )]
[( 1 0 1 )( 0 0 0 )]
[( 0 0 0 )( 1 0 1 )]
\end{verbatim}
\end{example}
In the {\tt Necklaces} package, the {\tt mkT} function  is  employed to recover 
the set $C(m,n)$ of all counting vectors ${\bm i}$ of length $m$ summing up to $n.$ Recall that the counting vector of an integer $n$ is a sequence $\bm i= (i_1, i_2, \ldots,i_m)$ of non-negative integers, named parts of $\bm i,$ such that
$i_1 + i_2 + \cdots + i_m= n,$ with $m \leq n.$ The integer $m=|{\bm i}|$ is  the  number of parts and is usually  called the length of the counting vector ${\bm i}$.
\begin{example}\label{2.3} {\rm To get all the counting vectors of $n=4$ in $m=2$ parts run}
\begin{verbatim}
> mkT(4,2, TRUE)
[( 2 )( 2 )]
[( 1 )( 3 )]
[( 3 )( 1 )]
[( 4 )( 0 )]
[( 0 )( 4 )]
\end{verbatim}
\end{example}
The {\tt Necklaces} package uses the set $C(m,n).$ Indeed, fix an alphabet $X$ of size $m$ and suppose to denote with  
\begin{enumerate}
\item[{\it i)}] $N[\bm i]$  (resp.  $B[\bm i]$)  the set of all representatives  of $m$-ary necklaces (resp. bracelets) of a fixed type $\bm i;$ 
\item[{\it ii)}] $N_m(n)$  (resp. $B_m(n)$) the set of all representatives of $m$-ary necklaces $N[\bm i]$  (resp. bracelets $B[\bm i]$) of type ${\bm i}$
with $|{\bm i}|=n$.
\end{enumerate}
Then we have 
\begin{equation}
N_m(n) = \bigcup_{{\bm i} \in C(m,n)} N[\bm i] \quad \hbox{and} \quad B_m(n) = \bigcup_{{\bm i} \in  C(m,n)} B[\bm i]
\label{eqNM}
\end{equation}
where the union is intended to be disjointed. The {\tt Necklaces} package provides the sets $N_m(n)$ and $B_m(n).$
\begin{example}
{\rm In Example  \ref{ex1}, for $X=\{0,1\},$ we have
$$N[(3,1)]=\{0001\}, N[(1,3)]=\{0111\}, N[(2,2)]=\{0011, 0101\}$$
$$N[(0,4)]=\{1111\}, N[(4,0)]=\{0000\}.$$
Using \eqref{eqNM} and Example \ref{2.3}, it is strainghforward to get  
$$N_2(4) = \{0000, 0001, 0011, 0101, 0111, 1111\}.$$
}
\end{example}
%------------------------------------------------------------------------
\section{Package necklaces in use}\label{sec:3} 
%------------------------------------------------------------------------
Given a vector in input, the {\tt cNecklaces} (resp. {\tt cBracelets} ) function generates the  necklace (resp. bracelet) whose this vector belongs to. When equal to {\tt TRUE}, the input flag {\tt bOut} permits to enumerate the vectors in output, and to print the elements of the class into a more compact form.  In addition,
the cardinality of the class is provided by the integer in parenthesis on the right of the last vector in the output list.
\begin{example} {\rm To print the elements of the necklace whose $001101$ belongs to, run}
\begin{verbatim}
> cNecklaces(c(0,0,1,1,0,1),TRUE)
[ 0 0 1 1 0 1 ]  ( 1 )
[ 0 1 0 0 1 1 ]  ( 2 )
[ 0 1 1 0 1 0 ]  ( 3 )
[ 1 0 0 1 1 0 ]  ( 4 )
[ 1 0 1 0 0 1 ]  ( 5 )
[ 1 1 0 1 0 0 ]  ( 6 )
\end{verbatim}
{\rm To print the elements of the bracelet whose $1021$ belongs to, run}
\begin{verbatim}
> cBracelets(c(1,0,2,1),TRUE)
[ 0 1 1 2 ]  ( 1 )
[ 0 2 1 1 ]  ( 2 )
[ 1 0 2 1 ]  ( 3 )
[ 1 1 0 2 ]  ( 4 )
[ 1 1 2 0 ]  ( 5 )
[ 1 2 0 1 ]  ( 6 )
[ 2 0 1 1 ]  ( 7 )
[ 2 1 1 0 ]  ( 8 )
\end{verbatim}
\end{example}
The set $N[\bm i]$ (resp. $B[\bm i]$) of all the representatives of necklaces (resp. bracelets) of a fixed type $\bm i$ (see Definition \ref{def1.1}) is generated running the function {\tt fNecklaces} (resp. {\tt fBracelets}).
\begin{example}
{\rm  To print the representatives of necklaces of type ${\bm i} =(2,1,1)$ run} 
\begin{verbatim}
> fNecklaces(c(2,1,1),TRUE)
[ 1 1 2 3 ]  ( 1 )
[ 1 1 3 2 ]  ( 2 )
[ 1 2 1 3 ]  ( 3 )
\end{verbatim} 
{\rm To print the representatives of bracelets of type ${\bm i} =(2,1,1)$ run} 
\begin{verbatim}
> fBracelets(c(2,1,1),TRUE)
[ 1 1 2 3 ]  ( 1 )
[ 1 2 1 3 ]  ( 2 )
\end{verbatim} 
{\rm Note that the elements of the class $\langle 1\,1\,2\,3 \rangle$ and $\langle 1\,1\,3\,2 \rangle$ (obtained from the first calling) are merged in just one bracelet  (obtained from the second calling) with representative  $1123.$  Indeed in the necklace
$ \langle 1\,1\,2\,3 \rangle=\{1123, 1231, 2311,$ $ 3112 \},$
the element  $1231\in \langle 1\,1\,2\,3 \rangle$ is obtained from  $1321 \in \langle 1\,1\,3\,2 \rangle =\{1132, 1321, 3211,$ $2113 \}$ by reversing the order and likewise for $2311,$ 
$3112, 1123 \in \langle 1\,1\,2\,3 \rangle$  from $1132, 2113, 3211\in \langle 1\,1\,3\,2 \rangle$ respectively. The input vector ${\bm i} =(2,1,1)$ refers to the alphabet $X=[3]$
and identifies $4$-lenght strings with $1$ occurring twice, $2$ and $3$ occurring once. To change the alphabet in $X=\{0,1,2\}$ - for example - add a last parameter equal to $0,$ as shown in the following example.}
\begin{verbatim}
> fNecklaces(c(2,1,1),TRUE,0)
[ 0 0 1 2 ]  ( 1 )
[ 0 0 2 1 ]  ( 2 )
[ 0 1 0 2 ]  ( 3 )
\end{verbatim}
\end{example}
Finally, the {\tt Necklaces} (resp. {\tt Bracelets}) function  returns the list of representatives of all necklaces (resp. bracelets) of length $n$ over the alphabet $X=[m].$ The function returns a vector: the second component contains the list and the first component indicates its length, as the following example shows.  
\begin{example}
{\rm  To get the number of representatives of binary necklaces of length $4$, run}
\begin{verbatim} 
> v=Necklaces(4,2)
> v[[1]]
[1] 6
\end{verbatim}
{\rm To get the elements of $N_2(4)$ run}
\begin{verbatim}
> for (w in v[[2]]) cat("[",unlist(w),"]\n");
[ 1 1 1 1 ]
[ 1 1 1 2 ]
[ 1 1 2 2 ]
[ 1 2 1 2 ]
[ 1 2 2 2 ]
[ 2 2 2 2 ]
\end{verbatim}
{\rm In agreement with \eqref{eqNM}, the following strategy is implemented in the {\tt Necklaces} (resp. {\tt Bracelets}) function. First the {\tt mkT} function  generates all the counting vectors of the integer $4$ in $2$ parts, that are ${\bm i} \in \{(2,2), (1,3), (3,1), (0,4), (4,0)\}$ see Example \ref{2.3}.  Then the {\tt fNecklaces} (resp. {\tt fBracelets}) function is used aiming to find all the representatives of type ${\bm i}.$
Thus for example one has}
\begin{verbatim}
> fNecklaces(c(2,2),TRUE)
[ 1 1 2 2 ]  ( 1 )
[ 1 2 1 2 ]  ( 2 )
> fNecklaces(c(1,3),TRUE)
[ 1 2 2 2 ]  ( 1 )
\end{verbatim}
{\rm and so on.}
\end{example}
Note that the cardinality of the sets $N_m(n)$ and $B_m(n)$ can also be calculated by running the  {\tt necklace} and {\tt bracelet} functions of the {\tt numbers} package  \cite{numbers}. 
\begin{example}\label{3.4}
{\rm  To compute the number of binary necklaces and bracelets of lenght $6$ given in Example \ref{ex1}, run}
\begin{verbatim} 
> Necklaces(6,2)[[1]]
[1] 14
> Bracelets(6,2)[[1]]
[1] 13
> library(numbers)
> necklace(2, 6)
[1] 14
> bracelet(2, 6)
[1] 13
\end{verbatim} 
\end{example}
%-----------------------------------------------------------------------------
\section{Applications}
%-----------------------------------------------------------------------------
 In most cases, a necklace has full size, which is equal to the length $n$ of a string. But, there are necklaces of size less than $n.$  In Example \ref{ex1}, the necklace $\{1000,0100,0010,0001\}$ has full size, but $\{1010,0101\}$ has size $2.$ In such a case the necklace is said periodic as its representative is a periodic string.  Recall that a string $\alpha$ is periodic if the string consists in a repetition of a single substring $\beta.$  Otherwise, the string is aperiodic. For example $0000$ is periodic of period $4,$ as 
$0000$ consists in a repeatition four times of the substring $0$. We might write $0000 = 0^4.$ Similarly $0101$ is periodic of period $2$ as  $0101 = (01)^2.$ In contrast, the string $0011 = (0011)^1$ is aperiodic. A necklace is said aperiodic if its representative is an aperiodic string. Next subsection \ref{LW} gives an example of construction of aperiodic necklaces using the {\tt necklaces} package. In subsection \ref{DB}, as a further possible application of the {\tt necklaces} package, de Bruijn sequences of order $n$ are generated. In such a case we use the aperiodic prefix of the representative of a necklace. Recall that the aperiodic prefix of a string $\alpha$ is its substring $\beta$ such that $\alpha= \beta^r$ for some positive integer $r.$ Hence for $0000$ the aperiodic prefix is $0,$ for $0011$ the 
aperiodic prefix is the same string $0011.$
\subsection{Lyndon words.}\label{LW}
{\rm A Lyndon word is an aperiodic necklace, that is a full size necklace with an aperiodic string as representative \cite{flajolet2009analytic}.  These words are involved in the expression of \eqref{genmom}. Additional applications involve compression analysis \cite{CA}, string matching \cite{crochemore1991two} and bioinformatics \cite{clare2019evolutionary}.

 To generate all Lyndon words of length $n$ over the alphabet $X=[m],$ run the {\tt cNecklaces} function and discard all the necklaces whose size is less than $n$. The {\tt R} code is:
\begin{verbatim} 
LyndonW <- function(n=1,m=2,bOut=FALSE, fn=1){
  if (n < 1 || m< 1) stop("n and m must be positive integers");
  i<-0; oL<-list(); iL<-Necklaces(n,m,fn);
  for (u in iL[[2]]) { 
    if (length(cNecklaces(u))==n) 
    {i<-i+1; oL[[i]]<-u};};
  oL<-lSort(oL);
  if (!bOut) return(oL);
  i<-1; for (w in oL) {cat("[",unlist(w),"]  (",i,")\n");i<-i+1};
}  
\end{verbatim} 
The {\tt LyndonW} function calls the {\tt lSort} function, taking in input a list of vectors and returning the same list ordered in lexicographical way. By default, the alphabet is $[m].$ Different alphabets $\{{\tt fn}, {\tt fn}+1, \ldots, {\tt fn}+m-1\}$ can be considered, changing the value of the input variable {\tt fn}.  In particular  if {\tt bOut = TRUE} the output is given in a compact form. 
\begin{example}
{\rm To generate all Lyndon words of length $3$ over the alphabet $[3]=\{1,2,3\},$ run}
\begin{verbatim} 
> LyndonW(3,3,TRUE)
[ 1 1 2 ]  ( 1 )
[ 1 1 3 ]  ( 2 )
[ 1 2 2 ]  ( 3 )
[ 1 2 3 ]  ( 4 )
[ 1 3 2 ]  ( 5 )
[ 1 3 3 ]  ( 6 )
[ 2 2 3 ]  ( 7 )
[ 2 3 3 ]  ( 8 )
\end{verbatim} 
\end{example}
%The same function can be used to generate all the bracelets of length $n$ over the alphabet $\{1,2,\ldots,m\},$  setting {\tt NB = "b"}.  
%\vskip0.01cm \noindent
%\textbf{\emph{Example 10:}} To generate all the bracelets of length $3$ over the alphabet $\{1,2,3\},$ run
%\begin{VerbIn} 
%>LyndonW(3,3,"b",TRUE)
%[ 1 1 2 ]  ( 1 )
%[ 1 1 3 ]  ( 2 )
%[ 1 2 2 ]  ( 3 )
%[ 1 2 3 ]  ( 4 )
%[ 1 3 3 ]  ( 5 )
%[ 2 2 3 ]  ( 6 )
%[ 2 3 3 ]  ( 7 )
%\end{VerbIn} 
%Note that the class {\tt [ 1 3 2 ]} is in the list of Example 9  but not in the list of Example 10, since the string $132$ is in the bracelet has length greater than $3.$
\subsection{de Bruijn sequences.} \label{DB}
de Bruijn sequences are cyclic $m$-ary strings in which every possible string of length $n$ of the alphabet occurs exactly once as a substring 
\cite{fredricksen1978necklaces}.  The following example clarifies the definition.
\begin{example}\label{exde}  {\rm For $n=4$ and $X=\{0,1\}$ ($m=2$), the minimum (in lexicographic order) de Bruijn sequence  is given in Fig. \ref{fig1}}.
\begin{figure}[h]
\includegraphics[width=3.3cm]{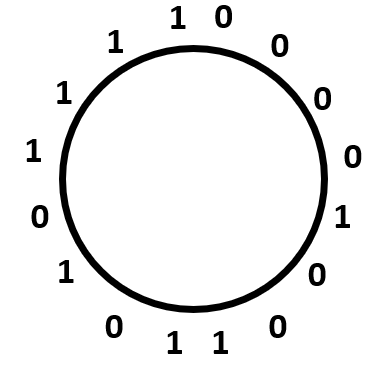}
\caption{A cyclic binary de Bruijn sequence of $4$-length strings.}
\label{fig1}
\end{figure}

 {\rm Table \ref{tab1} shows all the $2^4$ binary strings  (rows) of $\{0,1\}^4$ occuring (exactly once) as substring of the de Bruijn sequence of Fig. \ref{fig1}. In particular the columns of Table \ref{tab1} shows where each string  is located in the sequence: in bold the representatives of the binary necklaces of length $4.$ The de Bruijn sequence is given in the first row and the last $3$ digits have been added to reproduce the property of circularity.}

\begin{center}
\begin{table}[ht]
\begin{tabular}{ |cccccccccccccccc|ccc| } 
 \hline
 $0$&$0$&$0$&$0$&$1$&$0$&$0$&$1$&$1$&$0$&$1$&$0$&$1$&$1$&$1$&$1$&$0$&$0$&$0$\\ \hline
$\bf{0}$&$\bf{0}$&$\bf{0}$&$\bf{0}$& & & & & & & & & & & & & & &  \\
     &$\bf{0}$&$\bf{0}$&$\bf{0}$&$\bf{1}$& & & & & & & & & & & & & &  \\
     &     &$0$&$0$&$1$&$0$ & & & & & & & & & & & & &  \\
     &     &     &$0$&$1$&$0$&$0$& & & & & & & & & & & &  \\
     &     &     &     &$1$&$0$&$0$&$1$& & & & & & & & & & &  \\
     &     &     &     &     &$\bf{0}$&$\bf{0}$&$\bf{1}$&$\bf{1}$& & & & & & & & & &  \\
     &     &     &     &     &     &$0$&$1$&$1$&$0$& & & & & & & & &  \\
     &     &     &     &     &     &     &$1$&$1$&$0$&$1$& & & & & & & &  \\
    &     &     &     &     &     &      &     &$1$&$0$&$1$&$0$& & & & & & &  \\ 

    &     &     &     &     &     &      &     &     &$\bf{0}$&$\bf{1}$&$\bf{0}$&$\bf{1}$& & & & & &  \\ 
    &     &     &     &     &     &      &     &     &     &$1$&$0$&$1$&$1$& & & & &  \\ 
    &     &     &     &     &     &      &     &     &     &     &$\bf{0}$&$\bf{1}$&$\bf{1}$&$\bf{1}$& & & &  \\ 
    &     &     &     &     &     &      &     &     &     &     &     &$\bf{1}$&$\bf{1}$&$\bf{1}$&$\bf{1}$ & & &  \\ 
    &     &     &     &     &     &      &     &     &     &     &     &    &$1$&$1$&$1$&$0$ & &  \\ 
    &     &     &     &     &     &      &     &     &     &     &     &    &     &$1$&$1$&$0$&$0$ &  \\
    &     &     &     &     &     &      &     &     &     &     &     &    &     & &$1$&$0$&$0$&$0$  \\  
\hline
\end{tabular}
\caption{The occurrence of all $2^4$ strings of the set $\{0,1\}^4$ as substring of the de Bruijn sequence of Fig. \ref{fig1}.}
\label{tab1}
\end{table}
\end{center}
\end{example}
 Applications of de Bruijn sequences can be found in various fields. For example, computer programs make use of de Bruijn sequences when they play chess. A chess board is a square of $8 \times 8$ boxes that can be represented with integer numbers from $1$ to $64$  stored in a binary de Bruijn sequence of lenght $6$.  In chess programming \cite{leiserson1998using} the pieces of a given type are represented as a $64$-bit sequence and each bit flags the presence or absence of the piece type on a particular square of the chessboard \cite{frey1988creating}. In robotics these sequences are used to move the eye-camera of a robot. For example, if the space around the camera is labelled circularly as in Fig \ref{fig1}, then selecting a binary $4$-length string of Table \ref{tab1} a direction is identified in which moving the camera. Lastly, a further application is within modern DNA sequencing techniques. Usually the genome is broken down into small pieces (short reads) because it cannot be sequenced all at once. Available methods can only handle short DNA fragments from which the original sequence is reassembled, as if reconstructing a jigsaw puzzle. Various reassembly methods exist, among which De Bruijn sequences are widely used because they represent overlapping strings useful to manage overlapped short reads  \cite{compeau2011bruijn}.

There are various methods and tools to generate de Bruijn sequences (see for example \cite{ralston1982bruijn}). Necklaces are one of these tools relied on a concatenation of their aperiodic prefix as the following example shows.
\begin{example}\label{exde1}
{\rm Consider again Example \ref{exde} with $n=4, m=2$ and $X =\{0,1\}.$ To construct the de Bruijn sequence in Fig. \ref{fig1}, first generate the set $N_2(4)$ of all the representatives of binary necklaces of length $4,$ that is}
\begin{verbatim} 
> v=Necklaces(4,2,0)
> for (w in v[[2]]) cat("[",unlist(w),"]\n")
[ 0 0 0 0 ]
[ 0 0 0 1 ]
[ 0 0 1 1 ]
[ 0 1 0 1 ]
[ 0 1 1 1 ]
[ 1 1 1 1 ]
\end{verbatim} 
{\rm Then concatenate the aperiodic prefix of each representative (see Table \ref{tab2}) in the resulting sequence $0000100110101111$.}
\begin{center}
\begin{table}[ht]
\begin{tabular}{ |c|c| } \hline
{\rm Necklace} & {\rm Aperiodic prefix} \\ \hline
$\langle 0\,0\,0\,0 \rangle$ & $0$ \\
$\langle 0\,0\,0\,1 \rangle$ &  $0 0 0 1$ \\
$\langle 0\,0\,1\,1\rangle$ & $0 0 1 1$ \\
$\langle 0\,1\,0\,1\rangle$ & $0 1$ \\
$\langle 0\,1\,1\,1\rangle$ & $0 1 1 1$ \\
$\langle 1\,1\,1\,1 \rangle$ & $1$ \\ \hline 
\end{tabular}
\caption{In the first column the representatives  of all binary necklaces of length $4$ are given. In the second column their aperiodic prefix is given.}
\label{tab2}
\end{table}
\end{center}
\end{example}
The {\tt sBruijn} function generates the minimum string in lexicographic order using necklace concatenation \cite{MR3484729} as described in Example \ref{exde1}. The {\tt R} code is
\begin{verbatim} 
sBruijn <- function(n=1,m=2, fn=0, bSep=FALSE){
  if (n < 1) stop("n and m must be positive integers");
  sB<-"";sSep<-ifelse(bSep,".","");
  for (u in lSort(Necklaces(n,m,fn)[[2]])){                        
    apPrefix<-length(cNecklaces(u));
    sB<-paste0(sB,sSep,paste(u[1:apPrefix],collapse=""))}
  noquote(ifelse(bSep,substr(sB,2,nchar(sB)),sB))
}
\end{verbatim} 
In the code, first the function {\tt Necklaces} is ran aiming to generate all the representatives of $m$-ary necklaces of length $n.$ Then, after an increasing lexicographic ordering, the output is set by extracting their aperiodic prefix and concatenating them in a final sequence.  The alphabet is $\{0,1,\ldots,m-1\}$ by default.  Different alphabets such as $\{{\tt fn}, {\tt fn}+1, \ldots, {\tt fn}+m-1\}$ can be considered, changing the value of the input variable {\tt fn}.  
\begin{example}\label{44}{\rm To get the sequence in Fig. \ref{fig1} run}
\begin{verbatim} 
> sBruijn(4,2,0)
[1] 0000100110101111
\end{verbatim} 
{\rm Setting the input variable {\tt bSep=TRUE}, a separator is inserted among the output blocks to show the concatenation of the aperiodic prefixes.}
\begin{verbatim} 
> sBruijn(4,2,0,TRUE)
[1] 0.0001.0011.01.0111.1
\end{verbatim} 
\end{example} 
\begin{example}\label{inverse}
{\rm From the de Bruijn sequence, it is possible to generate the original overlapping strings as follows. Refer again to Table \ref{tab1}. Set $n=4$ to manage $4$-lenght strings and generate the binary de Bruijn sequence of Example \ref{44}. Then,
extract the first $n-1=3$ digits of the sequence, using the {\tt substr} function and concatenate these digits at the end of the sequence, after conversion to characters, with the {\tt paste0} function. In such a way,  the first row of Table \ref{tab1} is reproduced.  Lastly, 
repeatedly use  the {\tt substr} function to print all the $2^4$ substrings of length
$4.$ The {\tt R} code is the following:}
\begin{verbatim} 
> n=4; s<-sBruijn(n,2); s<-paste0(s,substr(s,1,n-1));
> for (i in 1:(nchar(s)-(n-1))) print(noquote(substr(s,i,i+(n-1))));
\end{verbatim} 
\end{example}
\begin{example}
{\rm Suppose $n=2, m=3$ and $X =[3].$ Then 
the set $N_3(2)$ of all the representatives of ternary necklaces of length $2$ is}
\begin{verbatim} 
> v=Necklaces(2,3)
> for (w in v[[2]]) cat("[",unlist(w),"]\n")
[ 1 1 ]
[ 1 2 ]
[ 1 3 ]
[ 2 2 ]
[ 2 3 ]
[ 3 3 ]
\end{verbatim} 
{\rm The concatenation consists in joining the aperiodic prefixes $1,12,13,2,23,3.$}
\begin{verbatim} 
> sBruijn(2,3,1,TRUE)
[1] 1.12.13.2.23.3
\end{verbatim} 
{\rm The final sequence is}
\begin{verbatim} 
> sBruijn(2,3,1)
[1] 112132233
\end{verbatim} 
{\rm Note that all possible strings of lenght $2$ are $11,12,21,13,32,22,23,33$  and they appear exactly once as substring of $112132233.$  To detect $31$, the sequence  $ 112132233$ should be closed circularly. Use the {\tt R} code of Example \ref{inverse} to recover all the $3^2$ strings of $\{1,2,3\}^2$ from the the de Bruijn sequence.}
\end{example}
%\vskip0.01cm \noindent
%\textbf{\emph{Example 11:}} To generate the minimum  de Bruijn sequences of order $3$ from cyclic binary strings, run 
%\begin{VerbIn} 
%> sBruijn(3,2)
%[1] 00010111
%> sBruijn(3,2,TRUE)
%[1] 0.001.011.1
%\end{VerbIn} 
%To change the alphabet, just run
%\begin{VerbIn} 
%> sBruijn(3,2,TRUE,1)
%[1] 1.112.122.2
%\end{VerbIn} 
%---------------------------------------------------------------------------------------------------
\section{Concluding remarks}
%---------------------------------------------------------------------------------------------------
The {\tt Necklaces} package was developed to list the elements of necklaces and bracelets under various hypothesis.  Further work would include the development of an algorithm in {\tt R}  to compute numerically  \eqref{genmom} and to list these sequences under more general hypothesis. The results in this note were obtained using {\tt R}~3.4.1. The 
{\tt kStatistics} package and the {\tt Necklaces} package are available from the Comprehensive {\tt R} Archive Network (CRAN) at \url{https://CRAN.R-project.org/}.

 \end{document}